\documentclass[1p]{elsarticle}
\usepackage{amssymb}
\usepackage{amsfonts}

\newtheorem{theorem}{Theorem}
\newtheorem{conjecture}[theorem]{Conjecture}

\newtheorem{lemma}[theorem]{Lemma}

\newproof{pf}{Proof}

\begin{document}
\title{Distant total irregularity strength of graphs 
via random vertex ordering}

\author{Jakub Przyby{\l}o\fnref{fn1,fn2}}
\ead{jakubprz@agh.edu.pl, phone: 048-12-617-46-38,  fax: 048-12-617-31-65}

\fntext[fn1]{Financed within the program of the Polish Minister of Science and Higher Education
named ``Iuventus Plus'' in years 2015-2017, project no. IP2014 038873.}
\fntext[fn2]{Partly supported by the Polish Ministry of Science and Higher Education.}

\address{AGH University of Science and Technology, al. A. Mickiewicza 30, 30-059 Krakow, Poland}

\begin{abstract}
Let $c:V\cup E\to\{1,2,\ldots,k\}$ be a (not necessarily proper) total colouring of
a graph $G=(V,E)$ with maximum degree $\Delta$.
Two vertices $u,v\in V$ are \emph{sum distinguished} if they differ with respect to sums of their incident colours,
i.e. $c(u)+\sum_{e\ni u}c(e)\neq c(v)+\sum_{e\ni v}c(e)$.
The least integer $k$ admitting such colouring $c$ under which
every $u,v\in V$ at distance $1\leq d(u,v)\leq r$ in $G$ are sum distinguished is denoted by ${\rm ts}_r(G)$.
Such graph invariants link the concept of the total vertex irregularity strength of graphs with
so called 1-2-Conjecture, whose concern is the case of $r=1$.
Within this paper we combine probabilistic approach with purely combinatorial one
in order to prove that ${\rm ts}_r(G)\leq (2+o(1))\Delta^{r-1}$ for every integer $r\geq 2$ and each graph $G$,
thus improving the previously best result: ${\rm ts}_r(G)\leq 3\Delta^{r-1}$.
\end{abstract}


\maketitle

\section{Introduction}
The cornerstone of the field of vertex distinguishing graph colourings is the graph invariant called \emph{irregularity strength}.
For a graph $G=(V,E)$ it is usually denoted by $s(G)$ and can be defined as the least integer $k$
so that we may construct an irregular multigraph,
i.e. a multigraph with pairwise distinct degrees of all vertices,
of $G$ by multiplying its edges, each at most $k$ times (including the original one), see~\cite{Chartrand}.
This study thus originated from the basic fact that no graph $G$ with more than one vertex is irregular itself
and related research on possible alternative definitions of an irregular graph, see e.g.~\cite{ChartrandErdosOellermann}.
Equivalently, $s(G)$ is also defined as the least $k$ so that there exists an edge colouring
$c:E\to\{1,2,\ldots,k\}$ such that for every pair $u,v\in E$, $u\neq v$,
the sum of colours incident with $u$ is distinct from the sum of colours incident with $v$.
Note that $s(G)$ exists only for graphs without isolated edges and with at most one isolated vertex.
It is known that $s(G)\leq n-1$, where $n=|V|$, for all such graphs, except for $K_3$, see \cite{Aigner,Nierhoff}.
This tight upper bound can however be improved in the case of
graphs with 
minimum degree $\delta\geq 1$ to $s(G)\leq 6\lceil\frac{n}{\delta}\rceil$ (what yields a better result whenever $\delta>12$ and for $\delta\in[7,12]$ if $n$ is larger than a small constant dependent on $\delta$), see
\cite{KalKarPf}, and to $s(G)\leq (4+o(1))\frac{n}{\delta}+4$ for graphs with
$\delta\geq n^{0.5}\ln n$, see \cite{MajerskiPrzybylo2}.
Many interesting results, concepts and open problems concerning this graph invariant
can also be found e.g. in~\cite{Bohman_Kravitz,Lazebnik,Dinitz,Faudree2,Faudree,Frieze,KalKarPf,Lehel,Przybylo},
and many others.

In~\cite{ref_BacJenMilRya1}, Ba\v ca et al. introduced a total version of the concept above.
Given any graph $G=(V,E)$ and a (not necessarily proper) total colouring $c:V\cup E\to\{1,2,\ldots,k\}$, let
\begin{equation}\label{weightofv}
w_c(v):=c(v)+\sum_{u\in N(v)}c(uv)
\end{equation}
denote the \emph{weight} of any vertex $v\in V$,
which shall also be called the \emph{sum at} $v$
and denoted simply by $w(v)$ in cases when $c$ is unambiguous from context.
The least $k$ for which there exists such colouring with $w(u)\neq w(v)$ for every $u,v\in V$, $u\neq v$,
is called the \emph{total vertex irregularity strength} of $G$ and denoted by ${\rm tvs}(G)$.
In~\cite{ref_BacJenMilRya1}, among others, it was proved that for every graph $G$
with $n$ vertices, $\lceil
\frac{n+\delta}{\Delta+1} \rceil \leq {\rm tvs}(G) \leq
n+\Delta-2\delta+1$.
Up to know the best upper bounds (for graphs with $\delta>3$) assert that ${\rm tvs}(G)\leq 3\lceil\frac{n}{\delta}\rceil+1$, see~\cite{AnhKalPrz},
and ${\rm tvs}(G)\leq (2+o(1))\frac{n}{\delta}+4$ for $\delta\geq n^{0.5}\ln n$, see~\cite{MajerskiPrzybylo1}.
Many other results e.g. for particular graph families can also be found in~\cite{AnholcerTotalForesets,NurdinEtAl,irreg_str2,RamdaniEtAl} and other papers.

In this article we consider a distant generalization of ${\rm tvs}(G)$ from~\cite{Przybylo_distant},
motivated among others by the study on distant chromatic numbers, see e.g.~\cite{DistChrSurvey} for a survey concerning these.
For any positive integer $r$, two distinct vertices at distance at most $r$ in $G$ shall be called \emph{$r$-neighbours}.
We denote by $N^r(v)$ the set of all $r$-neigbours of any $v\in V$ in $G$, and set $d^r(v)=|N^r(v)|$.
The least integer $k$ for which there exists a total colouring $c:V\cup E\to\{1,2,\ldots,k\}$
such that there are no $r$-neighbours $u,v$ in $G$ which \emph{are in conflict}, i.e. with $w(u)=w(v)$ (cf.~(\ref{weightofv})),
we call the \emph{$r$-distant total irregularity strength} of $G$, and denote by ${\rm ts}_r(G)$.
It is known that ${\rm ts}_r(G)\leq 3\Delta^{r-1}$ for every graph $G$
, see~\cite{Przybylo_distant},
also for a comment implying
that a general upper bound for ${\rm ts}_r(G)$ cannot be (much) smaller than $\Delta^{r-1}$.
In this paper we combine the probabilistic method with algorithmic approach similar to those in
e.g.~\cite{AnhKalPrz,Kalkowski12,MajerskiPrzybylo1,Przybylo_distant}
to prove that in fact ${\rm ts}_r(G)\leq (2+o(1))\Delta^{r-1}$ (for $r\geq 2$).

\begin{theorem}\label{Main_Theorem_tsr_probabil}
For every integer $r\geq 2$
there exists a constant $\Delta_0$ such that for each graph $G$ with maximum degree $\Delta\geq\Delta_0$,
\begin{equation}\label{main_total_distant_bound}
{\rm ts}_r(G)\leq 2\Delta^{r-1}+3\Delta^{r-\frac{4}{3}}\ln^2\Delta+4,
\end{equation}
hence
$${\rm ts}_r(G)\leq (2+o(1))\Delta^{r-1}$$
for all graphs.
\end{theorem}
It is also worth mentioning that the case of $r=1$ was introduced and considered separately in~\cite{12Conjecture},
where the well known 1-2-Conjecture concerning this invariant was introduced. It is known that ${\rm ts}_1(G)\leq 3$
for all graphs, see Theorem~2.8 in~\cite{Kalkowski12}, even in case of a natural list generalization of the problem, see~\cite{WongZhu23Choos},
though it is believed that the upper bound of $2$ should make the optimal general upper bound in
both cases,
 see~\cite{12Conjecture,PrzybyloWozniakChoos,WongZhuChoos}.

We also refer a reader to~\cite{Przybylo_distant_edge_probabil} to see an improvement of a similar probabilistic flavor for the upper bound from~\cite{Przybylo_distant} on the correspondent of ${\rm ts}_r(G)$ concerning the case of \emph{edge} colourings exclusively.

\section{Probabilistic Tools}
We shall use probabilistic approach in the first part of the proof of Theorem~\ref{Main_Theorem_tsr_probabil},
basing on
the Lov\'asz Local Lemma, see e.g.~\cite{AlonSpencer},
combined with the Chernoff Bound, see e.g.~\cite{JansonLuczakRucinski}
(Th. 2.1, page 26). We recall these below.
\begin{theorem}[\textbf{The Local Lemma}]
\label{LLL-symmetric}
Let $A_1,A_2,\ldots,A_n$ be events in an arbitrary pro\-ba\-bi\-li\-ty space.
Suppose that each event $A_i$ is mutually independent of a set of all the other
events $A_j$ but at most $D$, and that ${\rm \emph{\textbf{Pr}}}(A_i)\leq p$ for all $1\leq i \leq n$. If
$$ ep(D+1) \leq 1,$$
then $ {\rm \emph{\textbf{Pr}}}\left(\bigcap_{i=1}^n\overline{A_i}\right)>0$.
\end{theorem}
\begin{theorem}[\textbf{Chernoff Bound}]\label{ChernofBoundTh}
For any $0\leq t\leq np$,
$${\rm\emph{\textbf{Pr}}}({\rm BIN}(n,p)>np+t)<e^{-\frac{t^2}{3np}}~~{and}~~{\rm\emph{\textbf{Pr}}}({\rm BIN}(n,p)<np-t)<e^{-\frac{t^2}{2np}}\leq e^{-\frac{t^2}{3np}}$$
where ${\rm BIN}(n,p)$ is the sum of $n$ independent Bernoulli variables, each equal to $1$ with probability $p$ and $0$ otherwise.
\end{theorem}
Note that if $X$ is a random variable with binomial distribution ${\rm BIN}(n,p)$ where $n\leq k$, then we may still apply the Chernoff Bound above,
even if we do not know the exact value of $n$, to prove that
$\mathbf{Pr}(X>kp+t) < e^{-\frac{t^2}{3kp}}$ (for $t\leq\lfloor k\rfloor p$).

\section{Proof of Theorem~\ref{Main_Theorem_tsr_probabil}}
Fix an integer $r\geq 2$. Within our proof we shall not specify $\Delta_0$.
Instead, we shall assume that $G=(V,E)$ is a graph with sufficiently large maximum degree $\Delta$,
i.e. large enough so that all inequalities below are fulfilled.

We first partition $V$ into a subset of vertices with relatively small degrees and a subset of those with big degrees:
\begin{eqnarray}
S&=&\left\{u\in V:d(u)\leq \Delta^{\frac{2}{3}}\right\};\nonumber\\
B&=&\left\{u\in V:d(u)> \Delta^{\frac{2}{3}}\right\};\nonumber
\end{eqnarray}
Moreover, for every $v\in V$, we denote: $S(v)=N(v)\cap S$, $s(v)=|S(v)|$, $B(v)=N(v)\cap B$, $b(v)=|B(v)|$.

Now we randomly order the vertices of $V$ into a sequence.
For this goal,
associate with every vertex $v\in V$ a 
random variable $X_v\sim U[0,1]$ having the uniform distribution on $[0,1]$
where all these random variables $X_v$, $v\in V$ are independent,
or in other words pick a (real) number uniformly at random from the interval $[0,1]$ and associate it with $v$ 
for every $v\in V$.
Note that with probability one all these numbers are pairwise distinct.
In such a case, these
independent random variables uniquely define a natural ordering $v_1,v_2,\ldots,v_n$ of the vertices in $V$
where $X_{v_i}< X_{v_j}$ if and only if $1\leq i< j\leq n$.

For every vertex $v\in V$, any its neighbour or $r$-neighbour $u$ which precedes $v$ in the obtained
ordering of the elements of $V$ shall be called a \emph{backward neighbour} or \emph{$r$-neighbour}, resp., of $v$. Analogously, the remaining ones shall be called \emph{forward neighbours} or \emph{$r$-neighbours}, resp., of $v$, while the edges joining $v$ with its forward or backward neighbours shall be referred to as \emph{forward} or \emph{backward edges}, resp., as well.
Also, for any subset $W\subset V$, let $N_-(v)$, $N_-^r(v)$, $N_W^r(v)$ denote the sets of all backward neighbours, backward $r$-neighbours and $r$-neighbours in $W$ of $v$, respectively.
Set
$d_-^r(v)=|N_-^r(v)|$, $d_W^r(v)=|N_W^r(v)|$, and let $b_-(v)$ denote the number of backward neighbours of $v$ which belong to $B(v)$.

Denote $D(v)=\sum_{u\in N(v)}d(u)$
and note that
$d^r(v) \leq D(v)\Delta^{r-2} \leq d(v)\Delta^{r-1}$.

Let us also partition $V$ into a subset $I$ consisting of initial vertices of the obtained sequence and the remaining part $R$:
\begin{eqnarray}
I&=&\left\{v:X_v<\frac{\ln\Delta}{\Delta^{\frac{1}{3}}}\right\};\nonumber\\
R&=&\left\{v:X_v\geq\frac{\ln\Delta}{\Delta^{\frac{1}{3}}}\right\}.\nonumber
\end{eqnarray}

\begin{lemma}\label{MainSequencingLemma_total}
With positive probability, the obtained ordering has the following features for every vertex $v$ in $G$ with $b(v)\geq \Delta^{\frac{1}{3}}\ln\Delta$:
\begin{itemize}
\item[$F_1$:] $d^r_I(v) \leq 2d(v)\Delta^{r-\frac{4}{3}} \ln\Delta$;
\item[$F_2$:] if $v\in R$, then: $b_-(v)\geq X_v b(v)-\sqrt{X_v b(v)}\ln\Delta$;
\item[$F_3$:] if $v\in R$, then: $d^r_-(v)\leq X_v D(v) \Delta^{r-2}+\sqrt{X_vD(v)\Delta^{r-2}}\ln\Delta$.
\end{itemize}
\end{lemma}

\begin{pf}
For every vertex $v\in V$ of degree $d$ in $G$ and with $b(v)\geq \Delta^{\frac{1}{3}}\ln\Delta$ (hence also $d\geq \Delta^{\frac{1}{3}}\ln\Delta$), let $A_{v,1}$ denote the event that $d^r_I(v) > 2d\Delta^{r-\frac{4}{3}} \ln\Delta$,
let $A_{v,2}$ be the event that $v$ belongs to $R$ and $b_-(v) < X_v b(v)-\sqrt{X_v b(v)}\ln\Delta$, and
let $A_{v,3}$ denote the event that $v$ belongs to $R$ and $d^r_-(v) > X_v D(v) \Delta^{r-2}+\sqrt{X_vD(v)\Delta^{r-2}}\ln\Delta$.

As
$|N^r(v)| \leq d\Delta^{r-1}$ and for each $u\in N^r(v)$,
the probability that $u$ belongs to $I$ equals $\frac{\ln\Delta}{\Delta^{\frac{1}{3}}}$,
then by the Chernoff Bound (and the comment below it),
\begin{eqnarray}
\mathbf{Pr}(A_{v,1}) &\leq& \mathbf{Pr}\left(d^r_I(v) > d\Delta^{r-\frac{4}{3}} \ln\Delta + \sqrt{d\Delta^{r-\frac{4}{3}} \ln\Delta}\ln\Delta\right) \nonumber\\
&<& e^{-\frac{d\Delta^{r-\frac{4}{3}} \ln^3\Delta}{3d\Delta^{r-\frac{4}{3}} \ln\Delta}} = \Delta^{-\frac{\ln\Delta}{3}} < \frac{1}{\Delta^{3r}}.\label{Av1Ineq_total}
\end{eqnarray}

Subsequently note that for any $x\in [0,1]$:
\begin{eqnarray}
& & \mathbf{Pr}(b_-(v) < X_v b(v)-\sqrt{X_v b(v)}\ln\Delta | X_v=x)\nonumber\\
&=& \mathbf{Pr}({\rm BIN}(b(v),x) < x b(v)-\sqrt{x b(v)}\ln\Delta)\nonumber\\
&<& \frac{1}{\Delta^{3r}},\nonumber
\end{eqnarray}
where the last inequality follows by the Chernoff Bound if $\sqrt{xb(v)}\ln\Delta \leq x b(v)$, while it is trivial otherwise.
Hence,
\begin{equation}
\mathbf{Pr}(A_{v,2}) \leq \mathbf{Pr}(b_-(v) < X_vb(v)-\sqrt{X_v b(v)}\ln\Delta) \leq \int\limits_{0}^1\frac{1}{\Delta^{3r}}dx=\frac{1}{\Delta^{3r}}.\label{Av2Ineq_total}\end{equation}

For the sake of analyzing $A_{v,3}$, note now first that for $x\in[0, \frac{\ln\Delta}{\Delta^{\frac{1}{3}}})$,
\begin{equation}
\mathbf{Pr}(d^r_-(v) > X_v D(v) \Delta^{r-2}+\sqrt{X_vD(v)\Delta^{r-2}}\ln\Delta \wedge v\in R | X_v=x) = 0. \label{Av3Ineq_total_Part1}
\end{equation}
On the other hand, analogously as above, for $x\in [\frac{\ln\Delta}{\Delta^{\frac{1}{3}}},1]$:
\begin{eqnarray}
&&\mathbf{Pr}(d^r_-(v) > X_v D(v) \Delta^{r-2}+\sqrt{X_vD(v)\Delta^{r-2}}\ln\Delta \wedge v\in R | X_v=x) \nonumber\\
&\leq& \mathbf{Pr}({\rm BIN}(D(v) \Delta^{r-2},x) > x D(v) \Delta^{r-2}+\sqrt{xD(v)\Delta^{r-2}}\ln\Delta) \nonumber\\
&<&\frac{1}{\Delta^{3r}},\label{Av3Ineq_total_Part2}
\end{eqnarray}
where the last inequality follows by the Chernoff Bound, as $x\geq \frac{\ln\Delta}{\Delta^{\frac{1}{3}}}$ and $b(v)\geq \Delta^{\frac{1}{3}}\ln\Delta$ (where $D(v)\geq b(v)\Delta^{\frac{2}{3}}$) imply that $\sqrt{xD(v)\Delta^{r-2}}\ln\Delta \leq x D(v) \Delta^{r-2}$.
Hence, by~(\ref{Av3Ineq_total_Part1}) and~(\ref{Av3Ineq_total_Part2}),
\begin{equation}
\mathbf{Pr}(A_{v,3})
\leq \int\limits_{0}^1\frac{1}{\Delta^{3r}}dx=\frac{1}{\Delta^{3r}}.\label{Av3Ineq_total}\end{equation}

Note that each event $A_{v,i}$ is mutually independent of all other events except those $A_{u,j}$ with $u$ at distance at most $2r$ from $v$, $i,j\in\{1,2,3\}$,
i.e., at most $3\Delta^{2r}+2$ events. Thus, as by~(\ref{Av1Ineq_total}), (\ref{Av2Ineq_total}) and~(\ref{Av3Ineq_total}), the probability of each such event is bounded from above by $\Delta^{-3r}$,
by the Lov\'asz Local Lemma, with positive probability none of the events $A_{v,i}$ with $v\in V$ (and $b(v)\geq \Delta^{\frac{1}{3}}\ln\Delta$) and $i\in\{1,2,3\}$ appears.
\qed
\end{pf}

Let $v_1,v_2,\ldots,v_n$ be the ordering of the vertices of $V$ guaranteed by Lemma~\ref{MainSequencingLemma_total}.
Set
$$K=\Delta^{r-1} + \lceil\Delta^{r-\frac{4}{3}}\ln^2\Delta\rceil ~~~~{\rm and}~~~~ k=\lceil\Delta^{r-\frac{4}{3}}\ln^2\Delta\rceil,$$
and assign initial colour $1$ to all the vertices and initial colour $K+1$ to all the edges of $G$.
We shall construct our final colouring $f:V\cup E\to\{1,2,\ldots,2K+k+1\}$ using an algorithm within which
we shall be analyzing the consecutive vertices in the ordering (starting from $v_1$).
Denote by $c_t(a)$ the contemporary colour of every $a\in V\cup E$ at every stage of the ongoing algorithm
(hence initially $c_t(v)=1$ and $c_t(e)=K+1$ for every $v\in V$ and $e\in E$).
The final target sum of every vertex $v\in V$, $w_f(v)$, shall be chosen the moment $v$
is
 analyzed.
For every $v\in V$,
ever since $w_f(v)$ is chosen, we shall require so that
\begin{equation}\label{final-temporay_weight}
0\leq w_f(v)-w_{c_t}(v)\leq K.
\end{equation}
We shall admit at most two alterations of the colour for every edge in $E$ - only when any of its ends is being
analyzed
(vertex colours shall be adjusted at the end of the algorithm).
For every currently analyzed vertex
$v$ and its neighbour $u\in N(v)$,
we admit the following alterations of the colour of $e=uv$ (the moment $v$ is analyzed):
\begin{itemize}
\item adding $0,1,\ldots,K-1$ or $K$ if $e$ is a forward edge of $v$, $v\in S$ and $u\in B$,
\item adding $0,1,\ldots,k-1$ or $k$ if $e$ is a forward edge of $v$ (and $v\in B$ or $u\in S$),
\item adding $-K,-K+1,\ldots,K-1$ or $K$ if $e$ is a backward edge of $v$ and $u\in B$,
\item adding $-k,-k+1,\ldots,k-1$ or $k$ if $e$ is a backward edge of $v$ and $u\in S$,
\end{itemize}
so that afterwards (\ref{final-temporay_weight}) is fulfilled
for every vertex $u\in N_-(v)$ and for $v$ (after processing all edges incident with $v$).
Note that the admitted alterations guarantee that $c_t(e)\in\{1,2,\ldots,2K+k+1\}$ for every $e\in E$ at every stage of the construction.

Suppose we are about to analyze a vertex $v=v_i$, $i\in\{1,2,\ldots,n\}$, and thus far all our requirements have been fulfilled.
We shall show that in every case the admitted alterations on the edges incident with $v$ provide us more options for $w_{c_t}(v)$
than there are backward $r$-neighbours of $v$, and hence one of this options can be fixed as $w_f(v)$ so that this value is distinct from every $w_f(u)$ already fixed for any $u\in N^r_-(v)$. Denote the degree of $v$ by $d$,
and assume that $d>0$ (otherwise, we set $w_f(v)=1$):
\begin{itemize}
\item If $v\in I$, $v\in B$ and $b(v)\geq \Delta^{\frac{1}{3}}\ln\Delta$,
then the admitted alterations provide at least
$dk\geq d\Delta^{r-\frac{4}{3}}\ln^2\Delta$ available options for $w_{c_t}(v)$.
As by $F_1$ (from Lemma~\ref{MainSequencingLemma_total}), $|N^r_-(v)|\leq d_I^r(v)\leq 2d\Delta^{r-\frac{4}{3}} \ln\Delta < d\Delta^{r-\frac{4}{3}}\ln^2\Delta$,
at least one of these available options is distinct from all $w_f(u)$ with $u\in N^r_-(v)$.
\item If $v\in S$, then the admitted alterations provide at least
$s(v)k+b(v)K \geq s(v)\Delta^{r-\frac{4}{3}}\ln^2\Delta + b(v)(\Delta^{r-1}+\Delta^{r-\frac{4}{3}}\ln^2\Delta)$
available options for $w_{c_t}(v)$.
On the other hand,
$|N^r_-(v)|\leq d^r(v) \leq D(v)\Delta^{r-2} \leq (s(v)\Delta^{\frac{2}{3}}+b(v)\Delta)\Delta^{r-2}$,
hence at least one of these available options is distinct from all $w_f(u)$ with $u\in N^r_-(v)$.
\item If $v\in B$ and $b(v) < \Delta^{\frac{1}{3}}\ln\Delta$, then the admitted alterations provide at least
$dk\geq d\Delta^{r-\frac{4}{3}}\ln^2\Delta$ available options for $w_{c_t}(v)$.
On the other hand, analogously as in the case above,
$|N^r_-(v)|\leq d^r(v) \leq
s(v)\Delta^{\frac{2}{3}}\Delta^{r-2} + b(v)\Delta^{r-1} < d\Delta^{r-\frac{4}{3}}+\Delta^{r-\frac{2}{3}}\ln\Delta < d\Delta^{r-\frac{4}{3}}\ln^2\Delta$,
as $v\in B$ implies that $d\geq\Delta^{\frac{2}{3}}$.
We thus have at least one option available for $v$ distinct from all $w_f(u)$ with $u\in N^r_-(v)$.
\item If $v\in R$, $v\in B$ and $b(v) \geq \Delta^{\frac{1}{3}}\ln\Delta$, then by $F_2$ the number of available options for $w_{c_t}(v)$ via admitted alterations of colours of the edges incident with $v$ is not smaller than:
\begin{eqnarray}
b_-(v)K + (d-b_-(v))k &\geq& b_-(v)\Delta^{r-1} + d \Delta^{r-\frac{4}{3}}\ln^2\Delta \nonumber\\
&\geq& (X_v b(v)-\sqrt{X_v b(v)}\ln\Delta)\Delta^{r-1} + d \Delta^{r-\frac{4}{3}}\ln^2\Delta \nonumber\\
&\geq& X_v b(v)\Delta^{r-1} - \sqrt{d}\Delta^{r-1}\ln\Delta + d \Delta^{r-\frac{4}{3}}\ln^2\Delta \nonumber\\
&\geq& X_v b(v)\Delta^{r-1} + d \Delta^{r-\frac{4}{3}}\ln^2\Delta - d \Delta^{r-\frac{4}{3}}\ln\Delta \nonumber
\end{eqnarray}
(where the last inequality follows by the fact that $d\geq\Delta^{\frac{2}{3}}$).
This number is however greater than the number of backward $r$-neighbours of $v$, as by $F_3$,
\begin{eqnarray}
|N^r_-(v)| &\leq& X_v D(v) \Delta^{r-2}+\sqrt{X_vD(v)\Delta^{r-2}}\ln\Delta \nonumber\\
&\leq& X_v (b(v)\Delta +
s(v)\Delta^{\frac{2}{3}}) \Delta^{r-2}+\sqrt{d\Delta^{r-1}}\ln\Delta \nonumber\\
&\leq& X_v b(v)\Delta^{r-1} + d \Delta^{r-\frac{4}{3}} + d \Delta^{r-\frac{4}{3}}\ln\Delta. \nonumber
\end{eqnarray}
\end{itemize}
Thus in all cases there is at least one available sum, say $w^*$, for $v$ which is distinct from all $w_f(u)$ with $u\in N^r_-(v)$.
We then set $w_f(v)=w^*$ and perform the admitted alterations on the edges incident with $v$ so that $w_{c_t}(v)=w^*$ afterwards.

By our construction, after analyzing $v_n$, all $w_f(v_i)$ are fixed for $i=1,\ldots,n$ so that $w_f(u)\neq w_f(v)$ whenever $u$ and $v$ are $r$-neighbours in $G$ and~(\ref{final-temporay_weight}) holds for every $v\in V$.
We then modify (if necessary) the colour of every vertex $v$ by adding to it the integer $w_f(v)-w_{c_t}(v)$,
completing the construction of the desired total colouring $f$ of $G$ (by setting $f(a)=c_t(a)$ for every $a\in V\cup E$ afterwards).
Note that $1\leq f(e)\leq 2K+k+1$ for every $e\in E$ and, by~(\ref{final-temporay_weight}), $1\leq f(v)\leq K+1$ for every $v\in V$, hence the thesis follows.
\qed

\section{Remarks}
We have put an effort to optimize the second order term from the upper bound in~(\ref{main_total_distant_bound}),
up to a constant and a power in the logarithmic factor,
which could still be slightly improved (at the cost of the clarity of presentation).
Nevertheless, some multiplicative poly-logarithmic (in $\Delta$) factor seems unavoidable in this term within our approach.

We conclude by posing a conjecture, which to our believes
expresses
a true asymptotically optimal upper bound for the investigated parameters.
\begin{conjecture}\label{przybylo_main_con_total_distant}
For every integer $r\geq 2$ and each
graph $G$ with maximum degree $\Delta$,
$${\rm ts}_r(G) \leq (1+o(1))\Delta^{r-1}.$$
\end{conjecture}

\end{document}